%% file: pscc2022_template.tex
\documentclass{IEEEtran4PSCC}

\ifCLASSINFOpdf
   \usepackage[pdftex]{graphicx}
\else
   \usepackage[dvips]{graphicx}
\fi

\usepackage[cmex10]{amsmath}
\usepackage{algorithm}
\usepackage{algpseudocode} 
\usepackage{xcolor}

\usepackage{nomencl}
\makenomenclature

\makeatletter
\let\old@ps@headings\ps@headings
\let\old@ps@IEEEtitlepagestyle\ps@IEEEtitlepagestyle
\def\psccfooter#1{%
    \def\ps@headings{%
        \old@ps@headings%
        \def\@oddfoot{\strut\hfill#1\hfill\strut}%
        \def\@evenfoot{\strut\hfill#1\hfill\strut}%
    }%
    \def\ps@IEEEtitlepagestyle{%
        \old@ps@IEEEtitlepagestyle%
        \def\@oddfoot{\strut\hfill#1\hfill\strut}%
        \def\@evenfoot{\strut\hfill#1\hfill\strut}%
    }%
    \ps@headings%
}
\makeatother

\psccfooter{%
        \parbox{\textwidth}{\hrulefill \\ \small{PREPRINT Version: Submitted for consideration} \hfill \begin{minipage}{0.2\textwidth}\centering \vspace*{4pt} \includegraphics[scale=0.06]{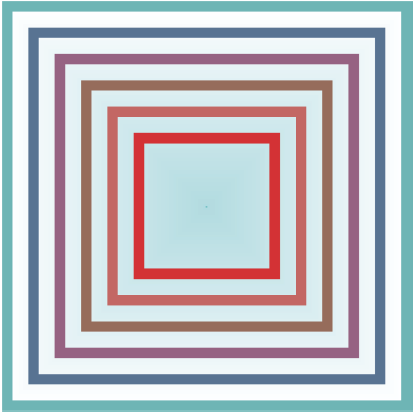}\\\small{PSCC 2022} \end{minipage} \hfill \small{Porto, Portugal --- June 27 -- July 1, 2022} \\ \small{to 22nd Power Systems Computation Conference}} %
}

\begin{document}
\newcommand{\ACOPFD}[0]{$\text{AC-OPF}^{\text{D}}$}
\newcommand{\ACOPFC}[0]{$\text{AC-OPF}^{\text{C}}$}
\newcommand{\ACOPFDs}[0]{$\text{AC-OPF}^{\text{D}}$ }
\newcommand{\ACOPFCs}[0]{$\text{AC-OPF}^{\text{C}}$ }

\newcommand{\xd}[0]{d}
\newcommand{\xdcont}[0]{d_c}
\newcommand{\xdbase}[0]{\hat{d}}
\newcommand{\xdtriv}[0]{\Tilde{d}}
\newcommand{\xdadj}[0]{\Delta d}
\newcommand{\xdadjcont}[0]{\Delta d_c}
\newcommand{\xdadjupper}[0]{\overline{\Delta d_c}}
\newcommand{\xdadjlower}[0]{\underline{\Delta d_c}}
\newcommand{\xdupper}[0]{\overline{d_c}}
\newcommand{\xdlower}[0]{\underline{d_c}}
\newcommand{\xext}[0]{x_c}

%
\title{Two-Stage Homotopy Method to Incorporate Discrete Control Variables into AC-OPF}

\author{
\IEEEauthorblockN{Timothy McNamara, Amritanshu Pandey, Aayushya Agarwal, Lawrence Pileggi}
\IEEEauthorblockA{Department of Electrical and Computer Engineering \\
Carnegie Mellon University\\
Pittsburgh, USA\\
\{tmcnama2, amritanp, aayushya, pileggi\}@andrew.cmu.edu}
}

\maketitle

\input{abstract}
\begin{IEEEkeywords}
Discrete controls, Robustness, Switched shunts, Tap changers, OPF
\end{IEEEkeywords}
\input{introduction}

\input{formulation}

\input{background}

\input{methodology}

\input{implementation}

\input{conclusion}

\input{ackowledgements}
\bibliographystyle{IEEEtran}
\bibliography{References}

\section{Appendix: Additional Results}
Table III contains a larger set of results obtained using the \ACOPFDs approach. For all these simulations a $k_{adj}$ value of 0.1 was used. The case file can be found at \cite{cases-github}. 
\textit{Note: this table was omitted from the original PSCC submission due to space constraints.}
\bigbreak

\begin{table}[h]
\onecolumn
\begin{center}
\caption{Additional \ACOPFDs results}
\begin{tabular}{|c|c|c|c|c|c|c|}
\hline
\textbf{Case} & \textbf{Buses} & \textbf{Tap Changers} & \textbf{Phase Shifter} & \textbf{Switched Shunt} & \textbf{Objective} & \textbf{\% Adj} \\ \hline
A  & 3022  & 981  & 4  & 399  & 5.377e5 & 24.9 \% \\ \hline
A* & 3022  & 981  & 4  & 399  & 5.357e5 & 25.6 \% \\ \hline
B  & 6867  & 759  & 5  & 161  & 1.216e5 & 90.8 \% \\ \hline
B* & 6867  & 759  & 5  & 161  & 1.216e5 & 93.9 \% \\ \hline
C  & 11152 & 530  & 0  & 500  & 5.294e5 & 77.8 \% \\ \hline
C* & 11152 & 530  & 0  & 500  & 5.294e5 & 82.2 \% \\ \hline
D  & 16789 & 997  & 2  & 1723 & 3.592e5 & 74.7 \% \\ \hline
D* & 16789 & 997  & 2  & 1723 & 3.591e5 & 89.2 \% \\ \hline
E  & 6549  & 1846 & 3  & 436  & 9.774e4 & 80.7 \% \\ \hline
E* & 6549  & 1846 & 3  & 436  & 9.777e4 & 91.1 \% \\ \hline
F  & 14393 & 0    & 5  & 724  & 2.317e4 & 33.7 \% \\ \hline
F* & 14393 & 0    & 5  & 724  & 2.312e4 & 33.2 \% \\ \hline
G  & 21849 & 997  & 2  & 1710 & 1.879e5 & 60.9 \% \\ \hline
G* & 21849 & 997  & 2  & 1710 & 1.879e5 & 66.4 \% \\ \hline
H  & 31156 & 12   & 36 & 2451 & 2.009e5 & 76.3 \% \\ \hline
H* & 31156 & 12   & 36 & 2451 & 2.009e5 & 63.9 \% \\ \hline
\end{tabular}
\end{center}
\end{table}

\end{document}

%% file: abstract.tex
\begin{abstract}
Alternating-Current Optimal Power Flow (AC-OPF) is an optimization problem critical for planning and operating the power grid. The problem is traditionally formulated using only continuous variables. Typically, control devices with discrete-valued settings, which provide valuable flexibility to the network and improve resilience, are omitted from AC-OPF formulations due to the difficulty of integrality constraints. We propose a two-stage homotopy algorithm to solve the AC-OPF problem with discrete-valued control settings. This method does not rely on prior knowledge of control settings or other initial conditions. The first stage relaxes the discrete settings to continuous variables and solves the optimization using a robust homotopy technique. Once the solution has been obtained using relaxed models, second homotopy problem gradually transforms the relaxed settings to their nearest feasible discrete values. We test the proposed algorithm on several large networks with switched shunts and adjustable transformers and show it can outperform a similar state-of-the-art solver.
\end{abstract}

%% file: introduction.tex
\section{Introduction}
\label{sec:introduction}

A critical framework for modeling and optimizing the efficiency of today’s power grid is based on the Alternating-Current Optimal Power Flow (AC-OPF) problem.
In the AC-OPF problem, a user-specified objective function, typically the cost of power generation, is optimized subject to network and device constraints.
These constraints include AC network constraints defined by Kirchhoff’s Voltage and Current Laws, as well inequality constraints representing operational limits, such as bounds on voltages, power generation, and power transfer, to ensure reliable operation of the power system.
It has been estimated that improved methods to model and run the US power grid could improve dispatch efficiency in the US electricity system leading to savings between \$6 billion and \$19 billion per year \cite{Cain}.
Moreover, improved AC-OPF solution techniques can improve the reliability and resiliency of the grid, which is under increasing duress from extreme weather events, such as California's ongoing wildfires and the aftermath of the winter 2021 storm in Texas.

Traditionally, the AC-OPF problem is a non-convex nonlinear problem with only continuous variables.
However, many devices deployed in the grid today have controls with discrete-valued settings, and are likely to become more widespread as the modernization of the grid continues.
Devices such as switched shunt banks and adjustable transformers can assist in balancing power flows in the system, meeting resilience-focused operational constraints, and locating a more optimal or resilient operating point than one located using fixed components settings.
The increased flexibility from these discrete devices also allows operators to avoid or delay costly upgrades to the network while increasing resiliency during extreme events. 
Given the significant potential benefits, a recent Grid Optimization (GO) Competition (Challenge 2), organized by ARPA-E, sought new robust approaches to AC-OPF where adjustments to discrete devices like tap changers, phase shifters, and switched shunt banks are included in the variable set \cite{go2}.

Although there are clear benefits to inclusion of discrete control devices in an AC-OPF study, doing so directly results in a mixed-integer non-linear program (MINLP) that is significantly harder to solve. 
Including these discrete settings in the variable space could lead to searching over a combinatorially-large solution space, which would be intractable to solve in practical time. 
If prior settings are at least known, then these can be used a starting point to begin the search for new settings.
However, prior settings may be not known in some situations like planning or policy studies. For instance, when engineers evaluate the feasibility of 50\% renewable penetration in a future U.S. Eastern Interconnection \cite{miso50}. 

One approach \cite{liu-linearization} to include these discrete variables relies on sequential linearization of the optimization problem and then handling the discrete variables using mixed-integer linear problem (MILP) techniques.
Unfortunately, the underlying network constraints can be highly non-linear with respect to certain control settings such as transformer tap ratios and phase shifts. Therefore these methods can suffer from a significant loss in model fidelity.
Additionally, linear relaxations to the OPF problem can lead to physically infeasible solutions, which are extremely undesirable for a grid dispatch \cite{Baker-DCOPF}.

The simplest approach to this obstacle is a two-stage rounding technique \cite{Tinney} \cite{Papalexopoulos}.
In this method, the discrete control settings are initially treated as continuous-valued, and a relaxed formulation of the AC-OPF problem is solved; then, these variables are fixed to their nearest respective discrete values, and the optimization is solved again with these variables held constant.
The first challenge with this approach is that for realistically-sized networks, solving the relaxed problem with continuous valued transformer taps, phase-shifters and switched shunts is computationally challenging, especially when good initial conditions are unavailable.
The second challenge is that the rounding step to map the continuous-valued settings to their nearest respective discrete values can result in a physically infeasible solution, which can be difficult to avoid when the discrete values are spaced far apart or when very many control variables are rounded at once. 
Directly rounding to the nearest discrete value also creates a discontinuous jump in the solution space, which is problematic for any Newton-method based solver that relies on first-order derivative continuity.

A number of methods have been introduced to address the two challenges above. 
Recent work  \cite{Coffrin} \cite{Lavei} has pushed the state-of-the-art for solving AC-OPF for realistic networks, but these formulations did not include discrete variables. 
For formulations that do include discrete control variables, several new approaches have been proposed \cite{liu-penalty-function,Macfie-discrete-shunt,Capitanescu,Murray-discrete} to eliminate non-physical solutions and degradation of optimality due to rounding.
\cite{liu-penalty-function} proposes utilizing a penalty term in the objective function to push each relaxed control variable towards an available discrete value, so that the disturbance introduced by the rounding step is smaller.
However, the use of penalty functions in optimizations with discrete variables can introduce stationary points and local minima \cite{discrete-opt-overview}.
Another approach is to select subsets of relaxed discrete variables to round in an outer loop, while repeatedly solving the AC-OPF problem in between subsets until all of the settings have been fixed.
In \cite{Macfie-discrete-shunt}, the authors present two methods for selecting which variables are rounded in each loop, and show these can reduce optimality degradation caused by rounding; however, unless only a single device is rounded at a time, this can introduce oscillations, since each rounding effectively adds a piecewise discontinuous function \cite{Katzenelson}.
\cite{Capitanescu} uses sensitivities with respect to the discrete variables as metric to determine when to round settings, with the help of either a merit function or a MILP solver.
In \cite{Murray-homotopy}, the authors point out that time constraints may make it impractical for grid operators to adjust a large number of control variable changes for a single dispatch.
To address this, they propose introducing a sparsity-inducing penalty term to the objective function along with a line-search of the discrete variable space to find a more limited number of control variable changes that can still improve optimality compared to holding all control variables constant.
However, this method inherently assumes knowledge of good settings, which might not be available in some use cases for AC-OPF, such as planning studies.

We propose a two-stage homotopy algorithm for solving AC-OPF problems that explains how to incorporate discrete controls variables, can scale to large networks, and is robust to potential lack of knowledge of prior settings. 
Our approach builds off of the homotopy-based AC-OPF methodology presented in \cite{Pandey-IMB}.
In the first stage,the discrete settings for the adjustable control devices are relaxed as continuous-valued variables and the optimization is solved using a robust homotopy technique.
The solution of the first stage is used to select discrete settings, and the respective variables are held constant thereafter. 
Then the errors induced by removing the relaxation are calculated and used in a second homotopy problem to locate a realistic solution.

The proposed novel approach can robustly determine a local optimum of any large real-world network with discrete controls without reliance on prior setting values. We also show that by choosing different homotopy paths, the proposed approach can obtain a variety of local minima solutions with significantly different generation dispatches.
In the results, we show that the method is not only novel in its approach but also more robust than another state-of-the-art optimization tool.


%% file: formulation.tex
\section{AC-OPF with Discrete Variables Formulation \ACOPFD}
\label{sec:formulation}
The paper solves the following AC-OPF with discrete control settings, which we will refer to as \ACOPFD:
\begin{subequations}
\begin{gather}
\underset{x,\xd}{\text{minimize   }}
f_0(x, \xd) \\
\text{subject to: }
 g(x, \xd) =0 \\
 h(x, \xd) \leq 0 \\
 \xd^i \in \mathcal{D}^i, \; i = 1, \ldots, n_d
\end{gather}
\end{subequations}
The vector $x$ consists of the continuous-valued variables of the optimization, including complex bus voltages, real and reactive power injections by generators, continuous shunts settings, and substitution variables to track transformer and line flows. 
The vector $\xd\in R^{n_d}$ represents the discrete control settings that are limited to a finite set of discrete values. In this paper, these include transformer taps $\tau$  and phase-shifters $\phi$ and discrete shunts $B^{sh}$ but the method is not restricted to just these. Each element of $\xd$, $\xd^i\in\xd$, has an associated integrality constraint (1d) restricting each device setting to be a finite value within the set, $\mathcal{D}^i$. The objective function (1a) is chosen depending on the purpose of the AC-OPF study, but typically represents the economic cost of power generation.
Constraint (1b) represents the AC network constraints, which can be formulated either to enforce net zero power-mismatch or current-mismatch at all nodes.
Constraint (1c) contains the bus voltage magnitude limits, branch and transformer thermal limits, and real and reactive power generation limits.

%% file: background.tex
\section{Background}
\label{sec:background}

\subsection{Homotopy Methods}
The Newton-Raphson (N-R) method is often used to solve the underlying non-linear equations in an AC-OPF formulation and its convergence can be sensitive to the choice of initial guess for the variables.
If the starting point is outside the basin of attraction for the problem, convergence can be very slow.
A class of successive-relaxation methodologies, known as homotopy methods, was introduced to mitigate such issues.
Homotopy is a numerical analysis method to solve a non-linear system of equations that traverses the solution space by deforming the non-linear equations from a trivial problem to the original.
The homotopy method initially defines a relaxation of the original problem which is trivial to solve, and proceeds to solve a sequence of deformations that ultimately leads back to the original problem.
Suppose $\mathcal{F}(x)=0$ is a set of non-linear equations that we aim to solve, we define a mapping to a trivial problem represented by $\mathcal{G}(x)=0$.
The deformation from the trivial problem to the original non-linear problem is controlled by embedding a  scalar homotopy factor $\nu \in [0,1]$ into the non-linear equations, thereby defining a  sequence of problems given by (\ref{eq:basic_homotopy}) .
\begin{equation}
\label{eq:basic_homotopy}
    \mathcal{H}(x, \nu) = \nu \mathcal{G}(x) + (1-\nu) \mathcal{F}(x)=0,~ \nu \in [0,1]
\end{equation}
Having determined the solution, $x^0$ to the trivial problem, $\mathcal{H}(x,1) = \mathcal{G}(x) = 0$,  we can iteratively decrease $\nu$, to move the system closer to the original problem, and use $x^0$ as our initial guess for the the next sub-problem.
By incrementally decreasing $\nu$ and solving the updated sub-problems, we traverse the solution space from the trivial problem to the original problem.
For this method to be effective, each sub-problem's solution should lie within the basin of attraction of the solution of the previous sub-problem, in order to exploit the quadratic convergence of N-R.
It is often challenging to develop a general homotopy method that ensures a proper traversal of the solution space, where there exists a feasible solution for every sub-problem $\mathcal{H}(x,\nu)=0$ as it traverses the path from $\nu: 1 \rightarrow 0 $ \cite{Allgower}. We present our homotopy method based on a circuit-inspired intuition that can intuitively ensure a feasible path. 

\subsection{Homotopy Methods in Power System}

A number of approaches have applied homotopy methods to solve power flow, AC-OPF, and other variants  in recent years \cite{Murray-homotopy,Pandey-IMB, Park-homotopy,Network-Stepping,Pandey-Tx-Stepping}.
In \cite{Pandey-Tx-Stepping}, the authors present a circuit-theoretic homotopy to robustly solve the power flow equations that embeds a homotopy factor in the equivalent circuits of the grid components to solve power flow.
These methods are extended in the Incremental Model Building (IMB) framework \cite{Pandey-IMB} to solve the AC-OPF optimization problem.
The idea behind the IMB framework is to build the grid from ground-up using an embedded homotopy factor in AC-OPF equations. 
IMB defines a relaxed problem, $\mathcal{G}(x)$, where the buses are almost completely shorted to one another, nearly all of the loads are removed at $\nu=1$, and $\nu$ is embedded in the generator limits so that the generator injections can be initially close to zero while remaining feasible with respect to the inequality constraints.
As a result, the relaxed network has very little current flow, so nearly all of the buses have voltage magnitude and angles close to that of the reference bus, and a flat-start initial point can reliably be used as a trivial solution of the homotopy sub-problem at $\nu=1$.
To satisfy the requirement that a feasible solution exist for every sub-problem along the homotopy path, fictitious slack current sources are introduced at each bus for sub-problems $\nu \neq 0$, and their injections are heavily penalized \cite{Jereminov-feasibility}.

While IMB shows an approach to solve AC-OPF without good initial conditions, it does not include discrete control variables in the formulation.
Even with a continuous relaxation to these variables, their introduction can significantly increase the nonlinearity of the network flow constraint equations and make the problem very challenging to solve.
In this work, we present a framework that builds on the IMB framework to include discrete control devices to solve the \ACOPFD robustly.

%% file: methodology.tex
\section{Two Stage Homotopy Method for Solving \ACOPFD}
\label{sec:methodology}
We propose a two-stage homotopy algorithm to solve the \ACOPFDs problem described in Section \ref{sec:formulation}.
The methods described are an extension of the IMB approach discussed above and we term the overall approach IMB+D.


The first stage, Stage I, applies a relaxation to the \ACOPFDs  problem in which we treat the discrete variables as continuous-valued by removing the integrality constraints. 
We refer to the resulting relaxed optimization as \ACOPFC, which is solved in Stage I.
We present a modeling framework for incorporating each adjustable device into the \ACOPFC~ homotopy formulation so as to preserve the IMB concept of slowly ``turning on'' the grid as the sequence of sub-problems is traversed.
After convergence of Stage I, in Stage II, the relaxed solution is used to select the nearest feasible discrete settings, and a second homotopy problem is defined to solve this problem. The local optima of this second stage yields an optimal solution for \ACOPFD with the discrete value constraints satisfied. 



\begin{figure}[h!]
    \centering
    \includegraphics[width=2.8in]{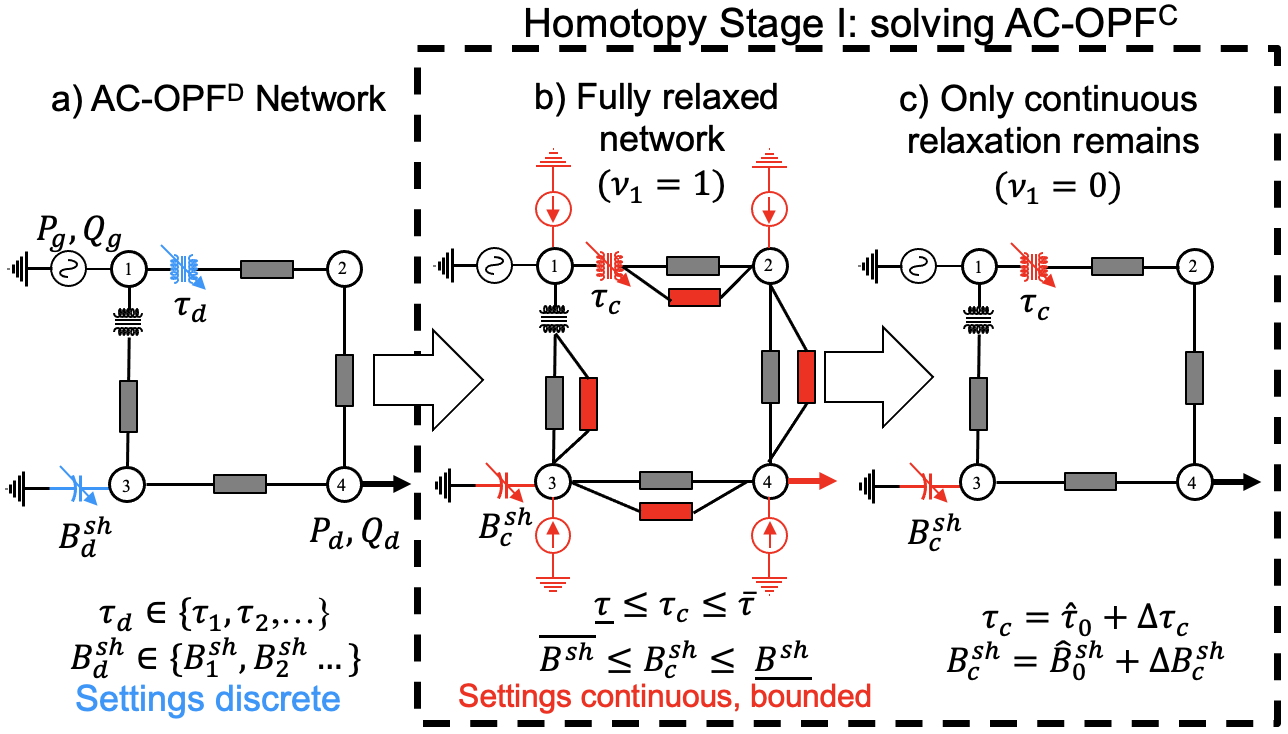}
    \caption{A simple network with discrete devices in blue (a). In the completely relaxed network (b), the relaxation elements are shown in red. At the end of Stage I (c), a solution is found with all relaxations removed except continuous settings.}
    \label{fig:stageI_figure}
\end{figure}
\subsection{Homotopy Stage I: Solving the Relaxed \ACOPFD}
\subsubsection{Embedding general network with homotopy}
In the constraints of an AC-OPF problem expressed using the current-voltage formulation, the current flows across lines are linear with respect to voltages, but the voltage and flow limit equations are quadratic, and the current injections from generators and loads are highly non-linear.
Even without the introduction of additional control devices, Newton method's may fail to converge if the initial guess for variables is not within the basin of attraction of a solution \cite{Murray-homotopy}. Therefore for large systems, without access to a reliable starting guess for Newton's method, we can ensure convergence using the IMB method presented in \cite{Pandey-IMB}.
By starting with a deformed version of the network in IMB that has very little current flow, there exists a high voltage solution in which the bus voltages are close to one another, which is close to a flat start guess.
To define a deformed network and have a smooth trajectory of intermediate deformations back to the original network, in Stage I the homotopy factor $\nu_1$ is embedded into the many of the parameters of the network's topology and devices: namely into high conductances in parallel to existing lines, the load factor, the generation limits (see Fig. \ref{fig:stageI_figure}.b). 

\subsubsection{Continuous relaxation and separation of discrete settings}
In Stage I, we apply a relaxation to the discrete control setting variables and split it into two components in order to leverage the robust methods of the IMB framework. 
First, a relaxation removing the integrality constraints (1d) is applied to replace the discrete variables $\xd$ with a continuous vector, $\xdcont$.
As a result, the constraints (1d) are replaced by:
\begin{equation}
    \xdlower^i \leq \xdcont^i \leq \xdupper^i \text{, } i \in 1...n_d
\end{equation}
where $\xdlower^i = \min(\mathcal{D}^i)$ and $\xdupper^i = \max(\mathcal{D}^i)$, representing the minimum and maximum possible settings for a device.

Even with this relaxation, adding devices that affect the flow of power across the network can significantly increase the nonlinearity of the constraint equations and the corresponding optimality conditions.
For example, introducing adjustable transformers causes the previously-linear transformer power flow model to be nonlinear with respect to voltage at its terminals, and adding phase shifters introduce trigonometric functions, and tap changers will introduce $\frac{1}{\tau}$ and $\frac{1}{\tau^2}$ terms.
To introduce these highly nonlinear models into the IMB framework, such that we preserve its initial trivial form where the entire grid is nearly shorted and a feasible homotopy path from $\nu_1 = 1 \rightarrow 0$, we design three measures.

First, each relaxed setting variable $\xdcont^i$ is separated into two components: a ``base'' value $\xdbase^i$ and a continuous-valued ``adjustment'' variable $\xdadjcont^i$:
\begin{equation}
\label{split}
    \xdcont^i = \xdtriv^i + \xdadjcont^i 
\end{equation}
Splitting $\xdcont^i$ affects its initialization and bounding, but the total value still drives the respective device behavior. 
However, this step allows maintaining the solution space of the variable $\xdadjcont$ around 0, which we have observed empirically improves the convergence in comparison to introducing $\xdcont$ as variable directly;  we believe this could be due to improved search directions for the Newton's method as the partials are dependent on the value of control variables.

Second, to maintain the trivial shorted form of the IMB method at $\nu_1 = 1$, we transform the discrete base setting $\xdtriv^i$ to be homotopy dependent such that during early stages of homotopy it has almost no impact on the network solution. For example, a tap ratio of 1.0 p.u. on a transformer would have no affect on the system as a whole. We define this smooth deformation of the discrete base setting $\xdtriv^i$ through embedding a homotopy factor $\nu_1$:
\begin{equation}
    \xdbase^i = \nu_1 \xdtriv^i + (1-\nu_1) \xdbase_0^i
\end{equation}
$\xdtriv^i$ is chosen as a setting value that would ensure that the trivial solution for the $\nu_1=1$ sub-problem is maintained. 
$\xdbase_0^i$ must be chosen from within the feasible domain. A good prior setting can be used, but the median value can always be used if the user does not know one.

Lastly, we ensure the effective settings do not stray from their respective trivial values in the early sub-problems of Stage I by adding a homotopy-dependent penalty term to the objective function parameterized by a $\nu_1$ and scaled by $k_{adj}$:
\begin{equation} 
    \label{modified_objective_kadj}
    f(x, \xdadjcont) = f_0(x, \xdadjcont) + \nu_1 k_{adj}  \sum_{i=1}^{N_d} |\xdadjcont^i|^2
\end{equation}
Including this term encourages minimizing $|\xdcont^i|$ more strongly in the early stages of homotopy. However, as $\nu_1$ is decreased, the penalty weight is reduced so that adjustment variables can move more freely as necessary to satisfy physical network constraints and decrease the primary objective function.
Note that all penalty terms have been removed from the objective function when the \ACOPFC ~is solved at $\nu_1=0$, so the form of the final \ACOPFC~ problem is independent of $k_{adj}$ value.


To address the likelihood of infeasible sub-problems along the homotopy path ($c(\nu_1),  \forall \nu_1 \in [0,1] $) where KCL cannot be satisfied without violating some variable limits, just as in IMB, homotopy dependent slack current injections (shown as red current sources in Fig. \ref{fig:stageI_figure}.b) are defined at each node in the network to allow satisfaction of conservation of charge \cite{Jereminov-feasibility}.
The magnitudes of the injection sources are penalized heavily in the objective function so that the sources only inject current if required for satisfaction of KCL, and the values are scaled by $\nu_1$ so that the fictitious sources are removed entirely when $\nu_1=0$.

\subsubsection{Solution of Stage I}
For each sub-problem defined in the homotopy path, the relaxed \ACOPFC~is solved using the primal-dual interior point (PDIP) approach \cite{Boyd}.
The Lagrangian for the sub-problem at any given $\nu_1 \in [0,1]$ is given by 
\begin{equation}
\begin{aligned}
    \label{eq:lagrangian}
    \mathcal{L}^{\nu_1}(\xext, \lambda,\mu) = f_0^{\nu_1}(\xext) + \lambda^T g^{\nu_1}(\xext) + \mu^T h^{\nu_1}(\xext)
\end{aligned}
\end{equation}
where $\xext = [x, \xdcont]^T$, $\lambda$ is vector of dual variables for the equality constraints, and $\mu$ is the vector of slack variables for the inequality constraints.
A local minimizer, $\theta^* = [\xext^*, \lambda^*, \mu^*]$, is sought by using Newton's method to solve for the set of perturbed first order KKT conditions:
\begin{equation}
\label{eq:KKT}
\begin{split}
&    \mathcal{F}(\theta)  = \begin{bmatrix}
\nabla_{\xext} f_0(\xext) + \nabla_{\xext}^T g(\xext) \lambda + \nabla_{\xext}^T h(x_{ext}) \mu \\
g(\xext) \\
\mu \odot h(\xext) + \epsilon
\end{bmatrix} = 0 \\
\end{split}
\end{equation}
where $\odot$ is element-wise multiplication.
In order to facilitate convergence and primal-dual feasibility of the solution, heuristics based on diode-limiting from circuit simulation methods are applied \cite{Pandey-Tx-Stepping}. 
The located $\theta^*$ is used as the initial guess for the next \ACOPFC sub-problem once the homotopy parameter (and thus the network relaxation) has been updated. 
This process is repeated until $\nu_1=0$, at which point a solution has been found to \ACOPFC (Fig. \ref{fig:stageI_figure}.c).

\subsection{Homotopy Stage II: Discretization}

After determining the optimal relaxed settings for the devices, we move onto the second stage of our approach, which solves for discrete value settings and the corresponding state of the grid using the relaxed solution from Stage I.
\subsubsection{Selection of discrete setting values}
A practical discrete setting value must be chosen for each control device.
For each control variable, the nearest-neighbor discrete value to  is selected. The chosen setting is evaluated to estimate whether snapping the variable to this value might result in infeasibility that prevents convergence.
The sensitivities of bus voltages with respect to the setting are calculated and an approximate voltage after rounding is calculated using the sensitivity vector and a first-order Taylor approximation.
With predicted voltage values, we check inequality constraints affected by the setting's perturbation.
If the chosen rounded value resulted in an infeasibility or violation of bounds, then second closest available setting is chosen and checked.

\begin{figure}[h]
    \centering
    \includegraphics[width=2.8in]{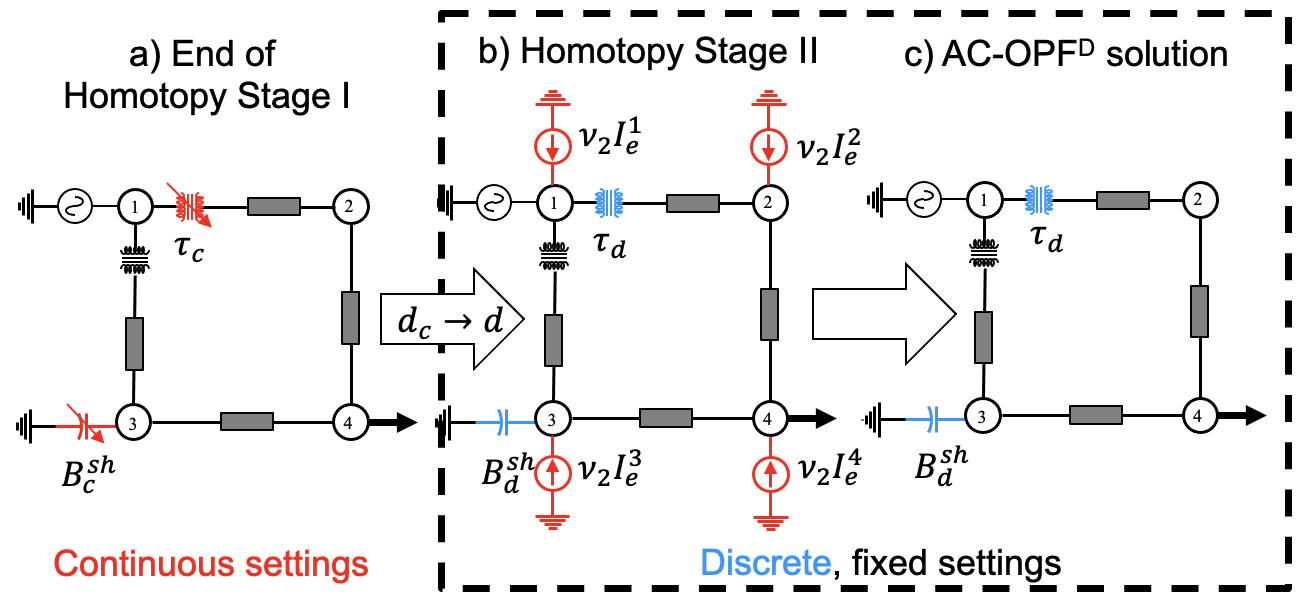}
    \caption{Process of using Stage I solution to discretize settings, formulate Stage II, and determine \ACOPFD solution}
    \label{fig:Stage-Two-Networks}
\end{figure}
\subsubsection{Stage II Homotopy: Error Injections}
At the termination of the Stage I, we have a solution vector $\theta^*$ that satisfies the perturbed KKT conditions for the continuous-valued OPF:
\begin{equation}
    \label{eq:KKT-sol-1}
    \mathcal{F}(\theta^{*}) = 0
\end{equation}
By changing the setting variables from their converged continuous values to realistic setting values ($\xdadjcont \rightarrow \xdadj$), the state vector is altered from $\theta^*$ to $\theta^\prime$.
Because of the adjustments, evaluating the KKT conditions at $\theta^\prime$ will result in violations of the conditions, which we will refer to as residual vector $R$:
\begin{equation}
    \label{eq:KKT-sol-res}
    \mathcal{F}(\theta^\prime) = R
\end{equation}
In the case of primal variables, we can think of R as a set of independent current sources that compensates for the current mismatch at each node due to rounding discrete device settings.
This idea can be extended to dual variables as well, as the underlying equations and nonlinearities have the same form \cite{Network-Stepping}.
Therefore, in the network disturbed by the rounding step, the power flow constraints can be satisfied immediately after rounding by adding a current sources to each of these perturbed buses, with current injection values defined by the current mismatches caused by rounding (see Fig. \ref{fig:Stage-Two-Networks}.b).
In a relaxed problem where we seek a solution to the same network but with the addition of $R$ to the respective equations, then we already know a solution to the problem: $\theta^\prime$.
\begin{equation}
    \mathcal{G}(\theta^\prime) = \mathcal{F}(\theta^\prime) - R = 0
\end{equation}
Therefore, we propose a second homotopy stage to find a feasible solution to the \ACOPFDs problem after control variables have been rounded, in which the relaxed system of equations are KKT conditions for the \ACOPFDs equations \eqref{eq:KKT} after rounding but with $R$ added as error injections.

\begin{equation}
\begin{split}
    \mathcal{H}_2(x,\nu_2) = (1-\nu_2)\mathcal{F}(\theta) + \nu_2(\mathcal{F}(\theta) - R) = 0 
\end{split}
\end{equation}
Here $\nu_2$, the Stage II scalar homotopy factor, is used to gradually reduce the residual injections, tracing a continuous path to a feasible solution for $\mathcal{F}(\theta)$ with rounded values.
This avoids taking a discontinuous jump between solving \ACOPFC and \ACOPFD.
When this homotopy problem is solved at $\nu_2=0$, the error injections have been removed, and we have located a realistic, feasible solution to the \ACOPFD.

\subsection{Generalization beyond IMB}
While the two-stage approach in this paper is described as an extension of the larger IMB framework \cite{Pandey-IMB} which assumes no knowledge of prior system settings, the two-stage homotopy algorithm can be applied to other approaches for solving AC-OPF without loss of generality. 
Consider the scenario in which good initial conditions for the general network are available but optimal settings for discrete variables ($\xd^*$) are unknown.
Such a use-case may occur if a grid planner wishes to re-evaluate existing settings for discrete devices (e.g. aims to move from a feasible setting $(\xd^k)$ to an optimal feasible setting $(\xd^*)$), or to explore the effects of upgrading of fixed devices to adjustable devices). 
In this situation, the grid planner may want to start from a feasible discrete setting but still explore whether a more optimal operating setting exists. In this scenario, we would still separate and relax each discrete setting $\xd^i$ according to \eqref{split}. 
But here, we account for this knowledge of good setting in Stage I of the algorithm by defining base values using the previously known setting $\xd^{k,i}$ such that $\xdbase^{i} = \xd^{k,i}$.

To ensure finding a solution at $\nu_1=1$ is simple, the objective function is modified according to \eqref{modified_objective_kadj}, but a very high penalty $k_{adj}$ value is used.
Now at $\nu_1=1$, with access to good initial conditions and feasible initial discrete settings, a trivial solution is obtained first.
However, as we traverse the homotopy path from $\nu_1 = 1 \rightarrow 0$, the adjustment values for discrete settings take  according to the objective $f_0(x, \xdadjcont)$ and eventually an optimal set of relaxed discrete settings are obtained at $\nu_1 = 0$. 
To find a feasible discrete setting, we perform Stage II as described in Section IV.B without modification. 

	    
		

%% file: implementation.tex
\section{Implementation and Evaluation}
\label{sec:implementation}
To test the effectiveness of the proposed algorithm, we run the \ACOPFDs solver on four networks based on cases used in the ARPA-E Grid Optimization Challenge 2 \cite{go2} that contain transformers with adjustable tap ratios or phase shifts, and discrete switched shunt banks.
The cases were modified to remove additional features of the GO formulation in order to focus on the efficacy of our approach for incorporating discrete control devices, and to make all generation costs linear.
The details of cases used are shown in Table I, and the files have been made available in a public Github repository \cite{cases-github}. 

\begin{table}[h]
\label{tab:case-info}
\caption{Properties of \ACOPFD ~Cases Tested}
\begin{tabular}{|c|c|c|c|c|c|}
\hline
\textbf{Name} & \textbf{Buses} & \textbf{Generators} & \textbf{Loads} & \textbf{Lines} & \textbf{Discrete Devices} \\ \hline
A & 3022 & 420 & 1625 & 3579 & 1384 \\ \hline
B & 6867 & 567 & 4618 & 7815 & 925 \\ \hline
C & 11152 & 1318 & 4661 & 16036 & 1030 \\ \hline
D & 16789 & 994 & 7846 & 23375 & 2722 \\ \hline
\end{tabular}
\end{table}
In each evaluation, $\xdbase$ values are set to the respective device settings listed in the .raw file, but to simulate running the \ACOPFDs without any prior knowledge of settings, a copy of each case file is created where each transformer's initial setting is listed as its median available setting, and each switched shunt bank has all switches off (0 p.u.).
These cases with initial settings removed are denoted with a * superscript.

\subsection{Robustness and Scalability of IMB+D}

Table \ref{tab:big-results} shows a summary of results.
$k_{adj}=0.1$ was applied equally across all devices in Stage 1, but normalized by the range of the individual device's settings ($\xdupper^i - \xdlower^i$).
As a comparison point, we also evaluated the same cases using the GravityX AC-OPF solver \cite{Gravity}, a leading submission to the GO Challenge which utilizes the IPOPT non-linear optimization tool \cite{IPOPT}.
``\% Adj'' indicates how many device settings in the solution differed from their initial value (the original .raw file values in the standard cases, and the simulated unknown in the * cases).
For each test, to evaluate the necessity of Stage II, a simple ``round and resolve'' approach is also attempted at the end of Stage I, in which the a feasible solution is sought immediately after rounding and fixing settings. 

The proposed approach is able to find solutions for all the cases, even when initial device settings are removed.
For Cases A and D, evaluating from the simulated unknown settings actually produces a very slightly better objective value.
GravityX produces a slightly better objective for both versions of cases B and C, but a less optimal solution for case A and does not converge for case D. 
Observe that Stage II is not strictly necessary for cases B, C, or D, but it is necessary for Case A to converge with discrete settings, which makes sense because in this case the tap steps are much further apart, so selecting discrete settings introduces a larger disturbance.
\begin{table}[]
\centering
\caption{Objective solutions and best $k_{adj}$ for tested cases}
\label{tab:big-results}
\begin{tabular}{|c|c|c|c|c|c|c|}
\hline
     & \multicolumn{3}{c|}{IMB+D ($k_{adj}=0.1$)}                   & \multicolumn{2}{c|}{GravityX+IP-OPT}                \\ \hline
Case  & Obj & \% Adj & Need Stage II? & Obj & \% Adj \\ \hline
A~ & 5.377e5 & 24.9\% & Yes & 5.561e5 & 23.8\% \\ \hline
A* & 5.357e5 & 25.6\% & Yes & 5.561e5 & 17.6\% \\ \hline
B~ & 1.216e5 & 90.8\% & No & 1.215e5 & 74.5\% \\ \hline
B* & 1.216e5 & 93.9\% & No & 1.215e5 & 71.7\% \\ \hline
C~ & 5.294e5 & 77.8 \% & No & 5.292e5 & 80.29\% \\ \hline
C* & 5.294e5 & 82.2\% & No & 5.292e5 & 69.4\% \\ \hline
D~ & 3.592e5 & 74.7\% & No & No Solution & N/A \\ \hline
D* & 3.591e5 & 89.2\% & No & No Solution  & N/A \\ \hline
\end{tabular}
\end{table}

\subsection{Impact of Different Homotopy Paths}
To construct different homotopy paths for the problem, we vary the Homotopy Stage I's adjustment penalty $k_{adj}$ parameter. We solve cases D and D* across a sweep of $k_{adj}$ values, and also at $k_{adj}=0$, meaning a penalty term is never used in Stage I.
Recall that the term is completely removed at $\nu_1 = 0$, so \textit{the same final problem is being solved}, regardless of choice of $k_{adj}$. 
Essentially, all that differs is the homotopy path to the original \ACOPFD problem. 

\begin{figure}[t]
    \centering
    \includegraphics[width=2.6in]{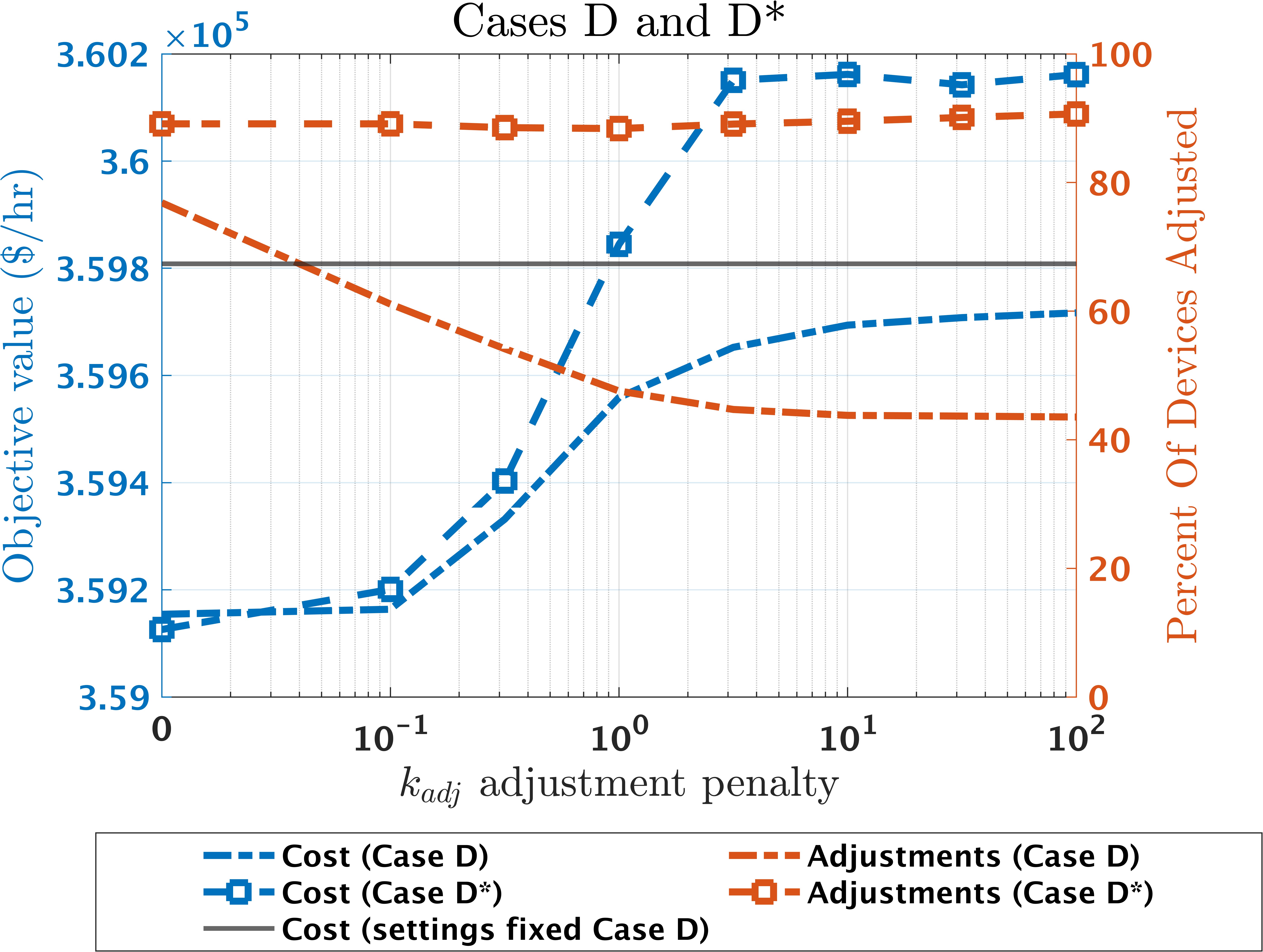}
    \caption{Changing $k_{adj}$ affects the homotopy path taken, and thus yields different solutions}
    \label{fig:stage1-graphs}
\end{figure}
Fig. \ref{fig:stage1-graphs} shows the effect of increasing $k_{adj}$ (i.e., varying the homotopy path) on both the final solution cost and percentage of adjustments made to devices in the solution.
This plot shows how it is possible to find multiple different local optima by parameterizing the penalty terms and essentially varying the homotopy path.
First, the we obtain the best objective function with very low penalty factors.
Additionally, when knowledge of base settings is available, increasing the $k_{adj}$ reduces the number of adjustments.
Surprisingly, however observe that for the D* simulations, increasing $k_{adj}$ does not have a large impact the number of discrete adjustments.
We hypothesize this is because some devices may require certain low or high settings for the network to be feasible, and so adjustments will be made regardless of the penalty value.

\subsection{Comparison of Local Solutions}

\begin{figure}
    \centering
    \includegraphics[width=2.5in]{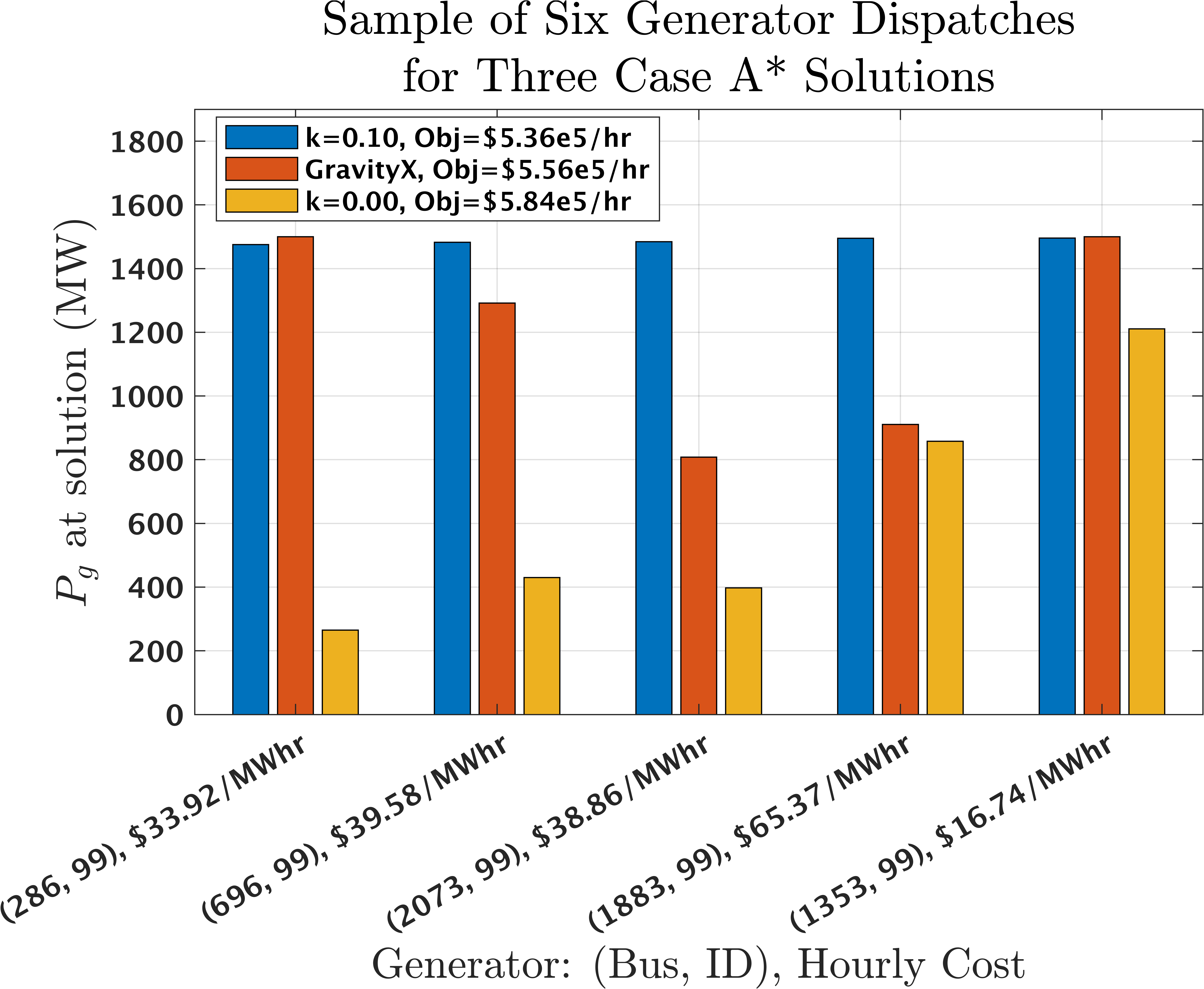}
    \caption{Changing the penalty scalar $k_{adj}$ on adjustments in Stage I causes a different homotopy path to be taken and can yield starkly different local optima in the final solution.}
    \label{fig:gen_dispatch}
\end{figure}

To further investigate the impact of different local solutions on the grid dispatch, we consider three sets of solution for the same exact problem as defined by Case A*: two solutions generated by IMB+D by taking separate homotopy paths through choice of $k_{adj}$, and one generated by GravityX.
The three solutions have generation dispatch costs in the range of \$5.36e5-\$5.84e5.
However, more interesting insights can be gathered by looking at the actual dispatch of six large generators in Case A* for these three distinct local optima, which are shown in Fig. \ref{fig:gen_dispatch}.
We notice that the dispatches for these generators vary significantly across the three dispatches.
This would imply that in the real world, the dispatch produced by an \ACOPFDs study can have widely varying patterns dependent on the method (e.g. Gravity-X vs. IMB+D) used to solve the non-convex problem, or even the choice of homotopy path within IMB+D based on the value of $k_{adj}$.
In real world, this could make it harder for grid operators to justify the choice of one dispatch over another as the global minima may not be easily obtainable.
Today, by running the convex DC-constraint based optimization with fixed settings, they are able to overcome the problem. 
However, due to the lack of granularity and accuracy, DC-based optimizations may be insufficient for future scenarios where adjustments to discrete devices are necessary in the optimization.

%% file: conclusion.tex
\section{Conclusion}
\label{sec:conclusion}
In this paper we developed a two-stage homotopy algorithm to robustly solve real-world AC-OPF problems while incorporating discrete control devices. The proposed approach uses fundamentals from circuit-theory to design the mechanisms of the underlying homotopy methods and the approach does not depend on access to good initial conditions to solve the overall problem. To evaluate this approach, we ran a series of tests for four networks containing transformers and shunts with discrete settings.
For two of the cases, this method performed better than a state-of-the-art solver. Furthermore, we showed that by constructing different homotopy paths, we find different local optima solutions for the same problem, which had significantly different generation dispatch patterns. Therefore, we believe that the nonconvexity of the solution space of \ACOPFDs problems warrant more investigation.



%% file: ackowledgements.tex
\section{Acknowledgements}
This work was sponsored in part by the  National Science Foundation under contract ECCS-1800812. 
This material is based upon work supported by the Google Cloud Research Credits program with the award GCP19980904.
We would like to sincerely thank Dr. Hassan Hijazi at Los Alamos National Laboratory for sharing a copy of his GravityX AC-OPF solver with us for comparison of our results. 